\documentclass{amsart}
\usepackage{amssymb}
\usepackage{graphicx}

\newtheorem{theorem}{Theorem}

\newtheorem{lem}{Lemma}

\newcommand{\A}{{\mathcal A}}
\newcommand{\U}{{\mathcal U}}

\newcommand{\PP}{{\mathcal P}}

\newcommand{\D}{{\mathbb D}}

\begin{document}
\bibliographystyle{amsplain}

\title[On a conjecture for the fifth coefficients for the class ${\U(\lambda)}$]{On a conjecture for the fifth coefficients for the class $\boldsymbol{\U(\lambda)}$}

\author[M. Obradovi\'{c}]{Milutin Obradovi\'{c}}
\address{Department of Mathematics,
Faculty of Civil Engineering, University of Belgrade,
Bulevar Kralja Aleksandra 73, 11000, Belgrade, Serbia}
\email{obrad@grf.bg.ac.rs}

\author[N. Tuneski]{Nikola Tuneski}
\address{Department of Mathematics and Informatics, Faculty of Mechanical Engineering, Ss. Cyril and Methodius
University in Skopje, Karpo\v{s} II b.b., 1000 Skopje, Republic of North Macedonia.}
\email{nikola.tuneski@mf.edu.mk}



\begin{abstract}
Let $f$ be function that is analytic in the unit disk $\D=\{z:|z|<1\}$, normalized such that $f(0)=f'(0)-1=0$, i.e., of type $f(z)=z+\sum_{n=2}^{\infty} a_n z^n$. If additionally,
\[ \left| \left(\frac{z}{f(z)}\right)^2 f'(z) -1\right|<\lambda \quad\quad  (z\in\D), \]
then $f$ belongs to the class $\U(\lambda)$, $0<\lambda\le1$.
In this paper we prove sharp upper bound of the modulus of the fifth coefficient of $f$ from $\U(\lambda)$ satisfying
\[ \frac{f(z)}{z}\prec \frac{1}{(1+z)(1+\lambda z)}, \]
("$\prec$" is the usual subordination) in the case when $0.400436\ldots \le\lambda\le1$.
\end{abstract}

\maketitle

\medskip

\section{Introduction and preliminaries}

\medskip

Let $\mathcal{A}$ consists of functions $f$ that analytic  in the open unit disc $\D=\{z:|z|<1\}$, with expansion
\begin{equation*}
f(z)=z+a_2z^2+a_3z^3+\cdots,
\end{equation*}
i.e., normalized such that $f(0)=f'(0)-1=0$. The famous Bieberbach conjecture from 1914 states that $|a_n|\le n$, $n=2,3,\ldots$, for the univalent functions from $\A$. The proof of the conjecture due to de Branges in 1985 \cite{Bra85} is one of the most celebrated results of the twentieth century.
Although, the conjecture is closed it remains an intriguing question to find upper bounds (preferably sharp) of the modulus of the coefficient for functons in various sublasses of univalent functions. One such class, that attracts significant attention in past decades is the class $\U(\lambda)$, $0<\lambda\le1$,
\[ \U(\lambda) = \left\{ f\in\A : \left| \left(\frac{z}{f(z)}\right)^2 f'(z) -1\right|<\lambda,\,z\in\D \right\}. \]
Functions from this class are proven to be univalent but not starlike which makes them interesting since the class of starlike functions is very wide. Overview of the most valuable results is given in Chapter 12 from \cite{DTV-book}.

\medskip

In \cite{OB-Monat},  the authors conjectured $|a_n|\le 1+\lambda+\lambda^2+\cdots+\lambda^{n-1}$ for the class $\U(\lambda)$ and $n\ge2$. In the same paper they proved that the conjecture is valid for $n=3$ and $n=4$, while for $n=2$ the proof is given in \cite{OB-Geom}. For the fifth coefficient the conjecture was proven in \cite{MONT-AMSJ-1} for the range  $2/3\le\lambda\le1$. The proofs for the third, fourth and the fifth coefficient rely on the claim from \cite{OB-Geom} that for every function $f$ from $\U(\lambda)$,
\begin{equation}\label{subord}
\frac{f(z)}{z}\prec \frac{1}{(1+z)(1+\lambda z)}.
\end{equation}
Here "$\prec$" denotes the usual subordination, i.e., $F(z)\prec G(z)$ for $F$ and $G$ analytic in $\D$, means that there exists function $\omega(z)$, also analytic in $\D$, such that  $\omega(0)=0$ and $|\omega(z)|<1$ for all $z\in \D$.

\medskip

Recently, in \cite{lipon}, by a counterexample, the authors showed that $f\in \U(\lambda)$ does not imply subordination \eqref{subord}. So, the cited estimates of $|a_n|$, for $n=3,4$ and 5, are correct only on the subclass of $\U(\lambda)$ consisting of functions satisfying the subordination \eqref{subord}.

\medskip

The estimate $|a_2|\le1+\lambda$ is correct and sharp on whole class $\U(\lambda)$ (see \cite{OB-Geom}).

\medskip

In this paper we study functions $f$ from $\U(\lambda)$ satisfying subordination \eqref{subord} and we extend the conjectured estimate for $n=5$ to the range $\lambda_0\le\lambda\le1$, where $\lambda_0=0.400436\ldots$ is the unique positive solution of the equation
\[9 \lambda^4-3 \lambda^3+\lambda^2+2 \lambda-1=0.\]

\medskip

For the proof we will use the following result for the class $\PP$ of Caratheodory functions, that are functions $p$ analytic in $\D$, of form $p(z)=1+p_1z+p_2z^2+\cdots$ with positive real part, i.e., $\operatorname{Re} p(z)>0$ for $z\in\D$. The result is due to Leverenz (\cite[Theorem 4(b)]{Lev}).

\medskip

\begin{lem}
Function $p(z)=1+p_1z+p_2z^2+\cdots$ has positive real part on the unit disk, if, and only if,
\begin{equation}\label{eq2}
\sum_{j=0}^\infty \left\{ \left| 2z_j+\sum_{k=1}^\infty p_kz_{k+j} \right|^2 - \left| \sum_{k=0}^\infty p_{k+1}z_{k+j} \right|^2 \right\} \ge 0
\end{equation}
for every sequence $\{z_k\}$ of complex numbers that satisfy $\lim_{k\to \infty} |z_k|^{1/k}<1$.
\end{lem}

\medskip

We will also need the following result by Prokhorov and Szynal \cite[Lemma 2, p.128]{prok}.

\begin{lem}\label{lem-prok}
Let $\omega(z)=c_1z+c_2z^2+c_3z^3+\cdots$ be analytic in $\D$ with $|\omega(z)|\le1$ for all $z\in\D$. If $\mu$ and $\nu$ are real numbers such that $2\le|\mu|\le4$ and $\nu\ge\frac{1}{12}(\mu^2+8)$, then $|c_3+\mu c_1c_2+\nu c_1^3|\le\nu$.
\end{lem}

\medskip

\section{Main result}

\medskip

If $p(z)=1+p_1z+p_2z^2+\cdots$ is a function from $\PP$, then there exists function $\omega(z)=c_1z+c_2z^2+\cdots$, analytic in $\D$, such that $\omega(0)=0$, $|\omega(z)|<1$ for all $z\in D$ and
\begin{equation*}
  p(z) = \frac{1+\omega(z)}{1-\omega(z)} \left( = 1+2\omega(z) +2\omega^2(z)+\cdots\right).
\end{equation*}
After comparing the coefficients we have
\begin{equation}\label{eq4}
\begin{split}
p_1 &= 2c_1,\\
p_2 &= 2(c_2+c_1^2),\\
p_3 &= 2(c_3+2c_{1}c_{2}+c_{1}^3),\\
p_4 &= 2(c_4+2c_{1}c_{3}+c_{2}^2+3c_1^2c_2+c_1^4).
\end{split}
\end{equation}

\medskip

This will help us to prove the main result.

\medskip
\begin{theorem}
Let $f(z)=z+a_2z^2+\cdots$ belongs to the class $\U(\lambda)$, where $\lambda_0\le\lambda\le1$, where $\lambda_0=0.400436\ldots$ is the unique positive solution of the equation
\[9 \lambda^4-3 \lambda^3+\lambda^2+2 \lambda-1=0.\]
If, additionally, $f$ satisfies subordination \eqref{subord}, then
\[ |a_5|\le 1+\lambda +\lambda^2+\lambda^3+\lambda^4,\]
and the result is sharp.
\end{theorem}

\medskip

\begin{proof}
If $f(z)=z+a_2z^2+\cdots\in \U(\lambda)$, $0<\lambda\le1$, satisfies subordination \eqref{subord}, then 
\begin{equation}\label{eq7}
\frac{f(z)}{z} = \frac{1}{(1-\omega(z))(1-\lambda\omega(z))} = 1+ \sum_{n=1}^\infty \frac{1-\lambda^{n+1}}{1-\lambda}\omega^n(z),
\end{equation}
where $\omega(z)=c_1z+c_2z^2+\cdots$ is analytic in $\D$, $|\omega(z)|<1$ for all $z\in\D$, and $\left.\frac{1-\lambda^{n+1}}{1-\lambda}\right|_{\lambda=1} = n+1$ for $n=1,2,\ldots$.
From \eqref{eq7} we have
\begin{equation}\label{eq-n7}
\begin{split}
a_5 &= (1+\lambda)c_4+2(1+\lambda+\lambda^2)c_1c_3+(1+\lambda+\lambda^2)c_2^2\\
&+3(1+\lambda+\lambda^2+\lambda^3)c_1^2c_2+(1+\lambda+\lambda^2+\lambda^3+\lambda^4)c_1^4.
\end{split}
\end{equation}

\medskip

On the other side, if we choose $z_k=0$ for $k>3$ in \eqref{eq2}, we have
\[
\begin{split}
& |2z_0+p_1z_1+p_2z_2+p_3z_3|^2 - |p_1z_0+p_2z_1+p_3z_2+p_4z_3|^2 \\
+& |2z_1+p_1z_2+p_2z_3|^2 - |p_1z_1+p_2z_2+p_3z_3|^2 \\
+& |2z_2+p_1z_3|^2 - |p_1z_2+p_2z_3|^2 +|2z_3|^2 - |p_1z_3|^2 \ge 0.
\end{split}
\]

\medskip

From here, we have that
\begin{equation}\label{eq-n8}
\begin{split}
 & L =: |p_1z_0+p_2z_1+p_3z_2+p_4z_3|^2 \\
  \le& R =: (|2z_0+p_1z_1+p_2z_2+p_3z_3|^2 - |p_1z_1+p_2z_2+p_3z_3|^2)\\
&\quad   + (|2z_1+p_1z_2+p_2z_3|^2 - |p_1z_2+p_2z_3|^2) \\
&\quad + (|2z_2+p_1z_3|^2  - |p_1z_3|^2) +|2z_3|^2.
\end{split}
\end{equation}
If we choose $p_1$, $p_2$, $p_3$ and $p_4$ from \eqref{eq4} and
\begin{equation*}\label{eq-9n}
\begin{split}
z_0 =& \lambda^2(1-\lambda)^2c_1^3,\\
z_1 =& \lambda^2c_2+(3\lambda^3-2\lambda^2)c_1^2,\\
z_2 =& 2\lambda^2c_1,\\
z_3 =& 1+\lambda,
\end{split}
\end{equation*}
then, after some calculations and comparing with \eqref{eq-n7}, for  $L$ defined in \eqref{eq-n8} we have that
\begin{equation}\label{eq-n10}
  L = 4|a_5|^2.
\end{equation}
Also, if we use that $|a+b|^2-|b|^2 = |a|^2+2\operatorname{Re}\{\overline{a}b\}$ ($a$ and $b$ are complex numbers), then by \eqref{eq-n8}:
\begin{equation}\label{eq-n11}
\begin{split}
R &= 4|z_0|^2 + 4\operatorname{Re}\{ (p_1z_1+p_2z_2+p_3z_3)\overline{z_0} \} +4|z_1|^2+ 4\operatorname{Re}\{(p_1z_2+p_2z_3)\overline{z_1}\}\\
& + 4|z_2|^2+4\operatorname{Re}\{(p_1z_3)\overline{z_2}\} + 4|z_3|^2.
\end{split}
\end{equation}
Considering each term of \eqref{eq-n11} by choosing the same values for $p_1$, $p_2$, $p_3$, $p_4$, $z_0$, $z_1$, $z_2$, $z_3$ as before, we receive:
\[
\begin{split}
|z_0|^2 &= \lambda^4(1-\lambda)^4|c_1|^6;\\
\operatorname{Re}\{& (p_1z_1+p_2z_2+p_3z_3)\overline{z_0} \} \le | p_1z_1+p_2z_2+p_3z_3|\cdot |\overline{z_0}| \\
 &= 2|(1+\lambda)c_3 + (3\lambda^2+2\lambda+2)c_1c_2+(3\lambda^3+\lambda+1)c_1^3| \lambda^2(1-\lambda)^2 |c_1|^3\\
 &= 2\lambda^2(1-\lambda)^2 |c_1|^3 (1+\lambda) \left| c_3+\left( 2+\frac{3\lambda^2}{1+\lambda} \right)c_1c_2+ \left( 1+\frac{3\lambda^3}{1+\lambda} \right)c_1^3  \right|\\
 &= 2\lambda^2(1-\lambda)^2 |c_1|^3 (1+\lambda) \left( 1+\frac{3\lambda^3}{1+\lambda} \right)  =  2\lambda^2(1-\lambda)^2(3\lambda^3+\lambda+1)|c_1|^3;
 \end{split}
\]
since in this case $2<|\mu|<4$ and $\nu\ge\frac{1}{12}(8+\mu^2)$ is equivalent with $\frac{\sqrt{52}-4}{9}=0.356789\ldots\le\lambda\le1$. Further,
\[
\begin{split}
|z_1|^2 &= |\lambda^2c_2+(3\lambda^3-2\lambda^2)c_1^2|^2 = \lambda^4|c_2|^2 + (3\lambda^3-2\lambda^2)^2|c_1|^4 \\
&\quad  + 2\operatorname{Re}\{\lambda^2(3\lambda^3-2\lambda^2)c_2\overline{c_1}^2\};\\
%
%
\operatorname{Re}\{&(p_1z_2+p_2z_3)\overline{z_1}\} \\
& = \operatorname{Re}\{ (4\lambda^2c_1^2+2(1+\lambda)(c_2+c_1^2))\cdot(\lambda^2\overline{c_2} +(3\lambda^3-2\lambda^2)\overline{c_1}^2) \} \\
& = \operatorname{Re}\{ (2(1+\lambda)c_2+2(2\lambda^2+\lambda+1)c_1^2)(\lambda^2\overline{c_2}+(3\lambda^3-2\lambda^2)\overline{c_1}^2) \} \\
& = \operatorname{Re}\{ 2\lambda^2(1+\lambda)|c_2|^2 + 2(2\lambda^2+\lambda+1)(3\lambda^3-2\lambda^2)|c_1|^4 \\
&\quad + (2\lambda^2(2\lambda^2+\lambda+1)+2(1+\lambda)(3\lambda^3-2\lambda^2))c_2\overline{c_1}^2 \};\\
|z_2|^2 &= 4\lambda^4|c_1|^2;\\
\operatorname{Re}\{&(p_1z_3)\overline{z_2}\}  = \operatorname{Re}\{2(1+\lambda)2\lambda^2|c_1|^2\} = 4\lambda^2(1+\lambda)|c_1|^2;\\
|z_3|^2 &= (1+\lambda)^2.
\end{split}
\]

\medskip

Using all previous facts and some transformations and calculations, from \eqref{eq-n11} we have
\[
\begin{split}
R &\le 4\left[ \left( \lambda^2(1-\lambda)^2|c_1|^3+3\lambda^3+\lambda+1 \right)^2 -(3\lambda^3+\lambda+1)^2 \right.\\
&\left.\quad +\lambda^2(\lambda^2+2\lambda+2)|c_2|^2 + (3\lambda^3-2\lambda^2)(3\lambda^3+2\lambda^2+2\lambda+2)|c_1|^4 \right.\\
&\left.\quad +4\lambda^2(\lambda^2+\lambda+1)|c_1|^2+2\lambda^2(3\lambda^3+3\lambda^2+2\lambda-1)\operatorname{Re}\{c_2\overline{c_1}^2\} +(1+\lambda)^2 \right].
\end{split}
\]
Since $3\lambda^3+3\lambda^2+2\lambda-1>0$ for $0.400436\ldots=\lambda_0\le\lambda\le1$, then $$2\lambda^2(3\lambda^3+3\lambda^2+2\lambda-1)\operatorname{Re}\{c_2\overline{c_1}^2\}\le 2\lambda^2(3\lambda^3+3\lambda^2+2\lambda-1)|c_2||c_1|^2,$$
and using that $|c_2|\le1-|c_1|^2$, we have
\[
\begin{split}
R &\le 4\left[ (\lambda^2(1-\lambda)^2|c_1|^3+3\lambda^3+\lambda+1)^2 +\lambda^2(\lambda^2+2\lambda+2)(1-|c_1|^2)^2 \right.\\
&\left. + (3\lambda^3-2\lambda^2)(3\lambda^3+2\lambda^2+2\lambda+2)|c_1|^4+4\lambda^2(\lambda^2+\lambda+1)|c_1|^2 \right.\\
&\left.+2\lambda^2(3\lambda^3+3\lambda^2 +2\lambda-1)(1-|c_1|^2)|c_1|^2+(1+\lambda)^2-(3\lambda^3+\lambda+1)^2 \right],
\end{split}
 \]
 and after some calculations, finally,
 \begin{equation*}\label{eq-12n}
   R \le 4 \left[ (\lambda^2(1-\lambda)^2|c_1|^3+3\lambda^3+\lambda+1)^2+F(\lambda,|c_1|^2) \right],
 \end{equation*}
 where
 \begin{equation}\label{eq-13n}
   F(\lambda,t) = 3\lambda^4(3\lambda^2-2\lambda-1)t^2+2\lambda^2(3\lambda^3+4\lambda^2+2\lambda-1)t-\lambda^2(9\lambda^4+5\lambda^2+4\lambda-2),
 \end{equation}
 $t=|c_1|^2$, $0\le t\le1$.

 \medskip

 If $\lambda=1$, then $F(1,t)=16(t-1)\le0$.

\medskip

If $0<\lambda<1$, then $3\lambda^2-2\lambda-1<0$ and the function $F(\lambda,t)$ attains its maximal value for
\[ t_0=\frac{\lambda^2(3\lambda^3+4\lambda^2+2\lambda-1)}{3\lambda^4(1+2\lambda-3\lambda^2)} \ge 1, \]
since this is equivalent to $9\lambda^4 -3\lambda^3 + \lambda^2+2\lambda-1\ge 0$, which is true because $\lambda_0\le\lambda\le1$. It means that $\max_{0\le t\le1} F(\lambda,1)=0$, i.e., $F(\lambda,|c_1|^2)\le0$ for all $\lambda_0\le\lambda \le1$ and $0\le|c_1|\le1$. By \eqref{eq-12n} we have
\begin{equation}\label{eq-14n}
\begin{split}
  R &\le 4\lambda^2(1-\lambda)^2|c_1|^3+3\lambda^3+\lambda+1)^2 \\
  &\le 4(\lambda^2(1-\lambda)^2+3\lambda^3+\lambda+1)^2\\
  &= 4(\lambda^4+\lambda^3+\lambda^2+\lambda+1)^2.
  \end{split}
\end{equation}
Finally, from \eqref{eq-n8}, \eqref{eq-n10} and \eqref{eq-14n} we have
\[ |a_5| \le 1+\lambda+\lambda^2+\lambda^3+\lambda^4.\]

\end{proof}

\end{document}